# A THERMODYNAMIC FRAMEWORK TO DEVELOP RATE-TYPE MODELS FOR FLUIDS WITHOUT INSTANTANEOUS ELASTICITY

SATISH KARRA AND K. R. RAJAGOPAL

*Dedicated to Prof. J. D. Goddard*

ABSTRACT. In this paper, we apply the thermodynamic framework recently put into place by Rajagopal and co-workers, to develop rate-type models for viscoelastic fluids which do not possess instantaneous elasticity. To illustrate the capabilities of such models we make a specific choice for the specific Helmholtz potential and the rate of dissipation and consider the creep and stress relaxation response associated with the model. Given specific forms for the Helmholtz potential and the rate of dissipation, the rate of dissipation is maximized with the constraint that the difference between the stress power and the rate of change of Helmholtz potential is equal to the rate of dissipation and any other constraint that may be applicable such as incompressibility. We show that the model that is developed exhibits fluid-like characteristics and is incapable of instantaneous elastic response. It also includes Maxwell-like and Kelvin-Voigt-like viscoelastic materials (when certain material moduli take special values).

## 1. Introduction

Recently, a thermodynamic framework has been put into place to describe the response of dissipative bodies that includes a large class of viscoelastic bodies (see [1] and [2], for details of the framework). With regard to the response of viscoelastic bodies, they consider the response to be that of a class of elastic bodies from an evolving set of configurations which they refer to as natural configurations (also see [3, 4]). The evolution of the natural configuration is determined by the rate of dissipation, or to be more precise, the maximization of the rate of dissipation. In a purely mechanical context, the response of the material is characterized by constitutively prescribing the stored energy (or Helmholtz potential) and the rate of dissipation functions. Since in a closed system the entropy increases to achieve its maximum equilibrium value, the quickest way in which the maximum could be reached is by maximizing the rate of dissipation. Of course, this is a plausible assumption and not a "principle". Also, to require such a criterion for an open system is not on very sound footing. However, it is surprising how well such a requirement works. Using such a thermodynamic framework a variety of material responses such as viscoelasticity [5, 6], twinning[7], solid-solid phase transition[8], plasticity[9], crystallization of polymers [10, 11], single crystal super alloys [12, 13], response of multi-network polymers [14] and anisotropic liquids [15] have been modelled. Particularly in viscoelasticity, this framework has been used to generalize one dimensional models due to Maxwell [16], Kelvin [17] and Voigt [18], Burgers [19], and the standard linear solid to three dimensions [5, 20, 21, 22]. Moreover, recently, it has been shown within the context of Maxwell fluid [23] that by choosing different forms for the stored energy and the rate of dissipation, one can obtain the same constitutive

*Key words and phrases.* viscoelastic fluid, instantaneous elasticity, Helmholtz potential, rate of dissipation, maximum entropy production, Lagrange multipliers .

KRR thanks the National Science Foundation for its support of his research.

relation for stress. It has also been shown (see [24]) that this framework leads to more than one three dimensional model that reduce within the context of one dimension to Burgers' fluid model.

Based on the aforementioned general framework, Rajagopal and Srinivasa [5] have developed a methodology to obtain rate-type models for fluids which possess instantaneous elasticity and in particular have derived a three dimensional generalization for the one dimensional Maxwell model (with the mechanical analog – a spring and a dashpot in series). When one assumes that the displacement gradients are small, this non-linear model leads to the classical upper convected Maxwell model. The aim of this paper is to use the framework to develop rate-type models for fluids which do not possess instantaneous elasticity. Specifically, we shall assume that the response from the natural configuration to the current configuration is akin to that of a generalized Kelvin-Voigt solid. Upon removal of external load, the body moves to the natural configuration with some "relaxation time" which is greater than the intrinsic time $t_m$ (also see [25]). By setting this relaxation time to a value less than $t_m$, one can obtain the class of models that can be generated using the framework in [5] and thus our work in this paper can be viewed as a generalization of the analysis in [5]. As an example, using our framework, we develop constitutive relations which in one dimension reduces to a dashpot and a Kelvin-Voigt element (which is a spring and a dashpot in parallel), in series (see Section 5). In carrying out the maximization of the rate of entropy production, one needs to decide what one maximizes with respect to. In most of the studies that have been cited, the maximization is with respect to appropriate kinematical variables and determining a constitututive representation for the stress. On the other hand, one could maximize with respect to the stress. In all the problems considered thus far, both methods lead to the same answer in explicit constitutive theories. However, it is possible to obtain a much larger class of constitutive relations following the latter procedure (see [26]). Also, from a philosophical perspective it is preferable to use the latter procedure. We shall not get into these issues here. In this paper we maximize with respect to the kinematical variables and obtain a constitutive expression for the stress.

While the idea of maximizing the rate of dissipation was also enunciated by Ziegler [27], Wehrli and Ziegler [28] our philosophy, interpretation and use of the maximization of the rate of entropy production is quite different, and the difference is not minor. We have discussed the differences in our approach in some detail in [26]. In fact, Ziegler's approach cannot be used to obtain a whole host of models including a whole class of implicit constitutive relations (see [26]). Also, the original work of Ziegler contains certain mathematical errors as pointed out in [2].

After a discussion of the preliminaries in the next section, we develop a general framework for the development of constitutive models for viscoelastic bodies that do not possess instantaneous elastic response in Section 3.1. Within our general framework, we derive a specific model in Section 3.2 which stores energy like a neo-Hookean solid with a rate of dissipation which depends on the stretching tensor of the natural configuration and the stretching tensor between the natural configuration and the current configuration. In Section 3.3, we show that our model reduces to either the Maxwell-like fluid or the Kelvin-Voigt-like solid under certain restrictions on the material parameters. We shall solve the problem of uniaxial extension in Section 4.1, followed by some remarks on the application of our framework to develop models for visco-elasto-plastic response in Section 6.

## 2. Preliminaries

Let $\kappa_R(\mathcal{B})$ and $\kappa_t(\mathcal{B})$ denote the reference configuration and the configuration of the body $\mathcal{B}$ at time $t$ (or the current configuration), respectively. Let $\mathbf{X}$ denote a typical point belonging to $\kappa_R(\mathcal{B})$ and $\mathbf{x}$ the same material point at time $t$, belonging to $\kappa_t(\mathcal{B})$. Let $\chi_{\kappa_R}$ denote a one to one mapping that assigns to each $\mathbf{X} \in \kappa_R(\mathcal{B})$, a point $\mathbf{x} \in \kappa_t(\mathcal{B})$, i.e.,

$$\mathbf{x} := \chi_{\kappa_R}(\mathbf{X}, t). \tag{2.1}$$

We shall assume that $\chi_{\kappa_R}$ is a sufficiently smooth mapping. The velocity $\mathbf{v}$, the velocity gradient $\mathbf{L}$ and the deformation gradient $\mathbf{F}_{\kappa_R}$ are defined through

$$\mathbf{v} := \frac{\partial \chi_{\kappa_R}}{\partial t}, \quad \mathbf{L} := \frac{\partial \mathbf{v}}{\partial \mathbf{x}}, \quad \mathbf{F}_{\kappa_R} := \frac{\partial \chi_{\kappa_R}}{\partial \mathbf{X}}. \tag{2.2}$$

It immediately follows that

$$\mathbf{L} = \dot{\mathbf{F}}_{\kappa_R} \mathbf{F}_{\kappa_R}^{-1}, \tag{2.3}$$

and the symmetric part of the velocity gradient $\mathbf{D}$ is given by

$$\mathbf{D} := \frac{1}{2}\left(\mathbf{L} + \mathbf{L}^T\right), \tag{2.4}$$

where $(.)^T$ denotes transpose of a second order tensor. The left and right Cauchy-Green stretch tensors $\mathbf{B}_{\kappa_R}$ and $\mathbf{C}_{\kappa_R}$ are defined through

$$\mathbf{B}_{\kappa_R} := \mathbf{F}_{\kappa_R} \mathbf{F}_{\kappa_R}^T, \quad \mathbf{C}_{\kappa_R} := \mathbf{F}_{\kappa_R}^T \mathbf{F}_{\kappa_R}. \tag{2.5}$$

Let $\kappa_{p(t)}$ denote the natural configuration associated with the configuration $\kappa_t$. We define $\mathbf{F}_{\kappa_{p(t)}}$ as the mapping from the tangent space at a material point in $\kappa_{p(t)}$ to the tangent space at the same material point at $\kappa_t$ (see Fig. 1). Similar to Eq. (2.5), we can also define the left and right Cauchy-Green tensors from the natural configuration to the current configuration [1]

$$\mathbf{B}_{\kappa_{p(t)}} := \mathbf{F}_{\kappa_{p(t)}} \mathbf{F}_{\kappa_{p(t)}}^T, \quad \mathbf{C}_{\kappa_{p(t)}} := \mathbf{F}_{\kappa_{p(t)}}^T \mathbf{F}_{\kappa_{p(t)}}. \tag{2.6}$$

The mapping $\mathbf{G}$ is defined through (see Fig. 1)

$$\mathbf{G} := \mathbf{F}_{\kappa_R \to \kappa_{p(t)}} := \mathbf{F}_{\kappa_{p(t)}}^{-1} \mathbf{F}_{\kappa_R}. \tag{2.7}$$

We define the tensors $\mathbf{B}_G$, $\mathbf{L}_G$ and $\mathbf{D}_G$ through

$$\mathbf{B}_G := \mathbf{B}_{\kappa_R \to \kappa_{p(t)}} = \mathbf{G}\mathbf{G}^T, \quad \mathbf{L}_G := \dot{\mathbf{G}}\mathbf{G}^{-1}, \quad \mathbf{D}_G := \frac{1}{2}\left(\mathbf{L}_G + \mathbf{L}_G^T\right). \tag{2.8}$$

In addition, let the tensors $\mathbf{L}_{p(t)}$ and $\mathbf{D}_{p(t)}$ be defined by

$$\mathbf{L}_{p(t)} := \dot{\mathbf{F}}_{\kappa_{p(t)}} \mathbf{F}_{\kappa_{p(t)}}^{-1}, \quad \mathbf{D}_{p(t)} := \frac{1}{2}\left(\mathbf{L}_{p(t)} + \mathbf{L}_{p(t)}^T\right). \tag{2.9}$$

Also, the principal invariants of $\mathbf{B}_{p(t)}$ are denoted by

$$I_{\mathbf{B}_{p(t)}} = tr(\mathbf{B}_{p(t)}), \quad II_{\mathbf{B}_{p(t)}} = \frac{1}{2}\left[(tr(\mathbf{B}_{p(t)}))^2 - tr(\mathbf{B}_{p(t)}^2)\right], \quad III_{\mathbf{B}_{p(t)}} = det(\mathbf{B}_{p(t)}), \tag{2.10}$$

where $tr(.)$ is the trace operator for a second order tensor and $det(.)$ is the determinant of a second order tensor.

---

[1] In this paper, henceforth we shall suppress $\kappa$ and denote $\mathbf{B}_{\kappa_{p(t)}}$ by $\mathbf{B}_{p(t)}$.

Now, from Eq. (2.7):

$$\dot{\mathbf{F}}_{\kappa_R} = \dot{\mathbf{F}}_{\kappa_{p(t)}}\mathbf{G} + \mathbf{F}_{\kappa_{p(t)}}\dot{\mathbf{G}}$$
$$\Rightarrow \dot{\mathbf{F}}_{\kappa_R}\mathbf{F}_{\kappa_R}^{-1} = \dot{\mathbf{F}}_{\kappa_{p(t)}}\mathbf{G}\mathbf{G}^{-1}\mathbf{F}_{\kappa_{p(t)}}^{-1} + \mathbf{F}_{\kappa_{p(t)}}\dot{\mathbf{G}} \quad (2.11)$$
$$\Rightarrow \mathbf{L} = \mathbf{L}_{p(t)} + \mathbf{F}_{\kappa_{p(t)}}\mathbf{L}_G\mathbf{F}_{\kappa_{p(t)}}^{-1}.$$

where $(\dot{.})$ is the material time derivative of the second order tensor given by

$$\dot{\mathbf{A}} := \frac{\partial(\mathbf{A})}{\partial t} + grad(\mathbf{A})[\mathbf{v}], \quad (2.12)$$

for a second order tensor $\mathbf{A}$ with $grad(.)$ being the gradient of a second order tensor with respect to the current configuration $\kappa_t$. Hence, from Eq. (2.11),

$$\mathbf{L}^T = \mathbf{L}_{p(t)}^T + \mathbf{F}_{\kappa_{p(t)}}^{-T}\mathbf{L}_G^T\mathbf{F}_{\kappa_{p(t)}}^T, \quad (2.13)$$

and so

$$\mathbf{D} = \mathbf{D}_{p(t)} + \frac{1}{2}\left[\mathbf{F}_{\kappa_{p(t)}}\mathbf{L}_G\mathbf{F}_{\kappa_{p(t)}}^{-1} + \mathbf{F}_{\kappa_{p(t)}}^{-T}\mathbf{L}_G^T\mathbf{F}_{\kappa_{p(t)}}^T\right]. \quad (2.14)$$

Further, we note that the upper convected Oldroyd derivative [29] (Also see [9] for the interpretation of Oldroyd derivative within the context of theory of multiple natural configurations) of $\mathbf{B}_{p(t)}$ can be related to $\mathbf{D}_G$ through (see [5])

$$\stackrel{\nabla}{\mathbf{B}}_{p(t)} := \dot{\mathbf{B}}_{p(t)} - \mathbf{L}\mathbf{B}_{p(t)} - \mathbf{B}_{p(t)}\mathbf{L}^T = -2\mathbf{F}_{\kappa_{p(t)}}\mathbf{D}_G\mathbf{F}_{\kappa_{p(t)}}^T. \quad (2.15)$$

Assuming that the body under consideration is incompressible, we shall record the balance of mass, and the balance of linear and angular momentum (in the absence of body couples):

$$div(\mathbf{v}) = 0, \quad \rho\dot{\mathbf{v}} = div(\mathbf{T}^T) + \rho\mathbf{b}, \quad \mathbf{T} = \mathbf{T}^T, \quad (2.16)$$

where $\rho$ is the density, $\mathbf{v}$ is the velocity, $\mathbf{T}$ is the Cauchy stress tensor, $\mathbf{b}$ is the specific body force and $div(.)$ is the divergence operator with respect to current configuration $\kappa_t$. In addition, the local form of balance of energy is

$$\rho\dot{\epsilon} = \mathbf{T}.\mathbf{L} - div(\mathbf{q}) + \rho r, \quad (2.17)$$

where $\epsilon$ denotes the specific internal energy, $\mathbf{q}$ denotes the heat flux vector and $r$ denotes the specific radiant heating.

Finally, we shall assume that the body under consideration undergoes isothermal processes. It is easy to modify the procedure to take into account non-isothermal processes. We shall invoke the second law of thermodynamics in the form of the *reduced energy dissipation equation*, given by (see [30, 2]):

$$\mathbf{T}.\mathbf{D} - \rho\dot{\psi} = \rho\theta\zeta := \xi \geq 0, \quad (2.18)$$

where $\mathbf{T}$ is the Cauchy stress, $\psi$ is the specific Helmholtz free energy, $\zeta$ is the rate of entropy production and $\xi$ is the rate of dissipation.[2]

---

[2]The terminology rate of dissipation usually refers to the conversion of working into energy in thermal form (heat). However, while considering general non-isothermal processes one uses the term to the product of the density, temperature and entropy production.

## 3. Constitutive assumptions

3.1. **General framework.** In this subsection, we shall first constitutively specify general forms for the specific Helmholtz potential $\psi$ and the rate of dissipation $\xi$. Using Eq. (2.18) and incompressbility as constraints, we shall maximize the rate of dissipation $\xi$, to obtain our constitutive relations. We shall assume that during its motion from $\kappa_{p(t)}$ to $\kappa_t$, the body stores energy as well as dissipates. The storage is due to elongation of the polymer networks. Now, assuming that the response of the body from $\kappa_{p(t)}$ to $\kappa_t$ is that of an isotropic body along with the assumption of incompressibility, we shall choose the specific Helmholtz free energy to be a function of the Cauchy-Green left stretch tensor $\mathbf{B}_{p(t)}$, i.e.,

$$\psi = \psi(\mathbf{B}_{p(t)}) = \hat{\psi}(I_{\mathbf{B}_{p(t)}}, II_{\mathbf{B}_{p(t)}}). \tag{3.1}$$

This implies that,

$$\dot{\psi} = 2\left[\left(\frac{\partial \hat{\psi}}{\partial I_{\mathbf{B}_{p(t)}}} + I_{\mathbf{B}_{p(t)}}\frac{\partial \hat{\psi}}{\partial II_{\mathbf{B}_{p(t)}}}\right)\mathbf{B}_{p(t)} - \frac{\partial \hat{\psi}}{\partial II_{\mathbf{B}_{p(t)}}}\mathbf{B}_{p(t)}^2\right].\mathbf{D}_{p(t)}. \tag{3.2}$$

Also, we shall assume the rate of dissipation to be a function of the stretching tensor and Cauchy-Green left stretch tensor between the $\kappa_{p(t)}$ to $\kappa_t$, the stretching tensor between $\kappa_R$ and $\kappa_{p(t)}$, i.e.,

$$\xi_m = \xi_m(\mathbf{B}_{p(t)}, \mathbf{D}_{p(t)}, \mathbf{D}_G). \tag{3.3}$$

In other words, the body dissipates energy during its motion from $\kappa_{p(t)}$ to $\kappa_t$ (due to continuous scission and healing of polymer networks) and also dissipates energy during its motion from $\kappa_R$ to $\kappa_{p(t)}$ (due to sliding of polymer strands over one another).

Due to the assumption of isotropic elastic response, the stored energy remains unchanged under any rotation. Hence, for our calculations, we shall assume that the natural configuration is rotated such that [3]

$$\mathbf{F}_{\kappa_{p(t)}} = \mathbf{V}_{\kappa_{p(t)}}, \tag{3.4}$$

and therefore

$$\begin{aligned}\mathbf{D}_{p(t)} &= \frac{1}{2}\left(\dot{\mathbf{V}}_{\kappa_{p(t)}}\mathbf{V}_{\kappa_{p(t)}}^{-1} + \mathbf{V}_{\kappa_{p(t)}}^{-1}\dot{\mathbf{V}}_{\kappa_{p(t)}}\right) \\ &= \frac{1}{2}\mathbf{V}_{\kappa_{p(t)}}^{-1}\dot{\mathbf{B}}_{p(t)}\mathbf{V}_{\kappa_{p(t)}}^{-1}.\end{aligned} \tag{3.5}$$

On substituting Eq. (3.2) into Eq. (2.18), we get

$$\mathbf{T}.\mathbf{D} - \mathbf{T}_{p(t)}.\mathbf{D}_{p(t)} = \xi_m(\mathbf{B}_{p(t)}, \mathbf{D}_{p(t)}, \mathbf{D}_G), \tag{3.6}$$

where

$$\mathbf{T}_{p(t)} := 2\rho\left[\left(\frac{\partial \hat{\psi}}{\partial I_{B_{p(t)}}} + I_{B_{p(t)}}\frac{\partial \hat{\psi}}{\partial II_{B_{p(t)}}}\right)\mathbf{B}_{p(t)} - \frac{\partial \hat{\psi}}{\partial II_{B_{p(t)}}}\mathbf{B}_{p(t)}^2\right]. \tag{3.7}$$

Substituting Eq. (2.14) into Eq. (3.6), we obtain that

$$\frac{1}{2}\mathbf{T}.\left(\mathbf{F}_{\kappa_{p(t)}}\mathbf{L}_G\mathbf{F}_{\kappa_{p(t)}}^{-1} + \mathbf{F}_{\kappa_{p(t)}}^{-T}\mathbf{L}_G^T\mathbf{F}_{\kappa_{p(t)}}^T\right) + \left(\mathbf{T} - \mathbf{T}_{p(t)}\right).\mathbf{D}_{p(t)} = \xi_m(\mathbf{B}_{p(t)}, \mathbf{D}_{p(t)}, \mathbf{D}_G). \tag{3.8}$$

Since $\mathbf{T}$ is symmetric, Eq. (3.8) reduces to

$$\mathbf{F}_{\kappa_{p(t)}}^T\mathbf{T}\mathbf{F}_{\kappa_{p(t)}}^{-T}.\mathbf{L}_G + \left(\mathbf{T} - \mathbf{T}_{p(t)}\right).\mathbf{D}_{p(t)} = \xi_m(\mathbf{B}_{p(t)}, \mathbf{D}_{p(t)}, \mathbf{D}_G). \tag{3.9}$$

---

[3] For the application of this thermodynamic framework to anisotropic fluids the reader should refer to [15].

Also the assumption of incompressibility leads to

$$tr(\mathbf{D}) = tr(\mathbf{D}_{p(t)}) = tr(\mathbf{D}_G) = 0. \tag{3.10}$$

Now, following Rajagopal and Srinivasa [5], we maximize the rate of dissipation $\xi_m$ by varying $\mathbf{L}_G, \mathbf{D}_{p(t)}$ for fixed $\mathbf{B}_{p(t)}$ with Eqs. (3.9), (3.10) as constraints. We maximize the auxillary function $\Phi$ given by

$$\Phi := \xi_m + \lambda_1 \left[\xi_m - \mathbf{F}^T_{\kappa_{p(t)}} \mathbf{T} \mathbf{F}^{-T}_{\kappa_{p(t)}} . \mathbf{L}_G - \left(\mathbf{T} - \mathbf{T}_{p(t)}\right) . \mathbf{D}_{p(t)}\right] \\ + \lambda_2 (\mathbf{I}.\mathbf{D}_G) + \lambda_3 (\mathbf{I}.\mathbf{D}_{p(t)}), \tag{3.11}$$

where $\lambda_1, \lambda_2, \lambda_3$ are Lagrange multipliers. By setting, $\partial\Phi/\partial\mathbf{D}_{p(t)} = 0$ and $\partial\Phi/\partial\mathbf{L}_G = 0$, we get

$$\mathbf{T} = \mathbf{T}_{p(t)} + \frac{\lambda_3}{\lambda_1}\mathbf{I} + \left(\frac{\lambda_1 + 1}{\lambda_1}\right) \frac{\partial \xi_m}{\partial \mathbf{D}_{p(t)}}, \tag{3.12}$$

and

$$\mathbf{T} = \frac{\lambda_2}{\lambda_1}\mathbf{I} + \left(\frac{\lambda_1 + 1}{\lambda_1}\right) \mathbf{F}^{-T}_{\kappa_{p(t)}} \frac{\partial \xi_m}{\partial \mathbf{L}_G} \mathbf{F}^T_{\kappa_{p(t)}}. \tag{3.13}$$

At this juncture it is worth recalling the comments in the introduction concerning maximization with respect to the kinematical quantities. On substituting Eqs. (3.12), (3.13) in Eq. (3.9), we finally obtain

$$\left(\frac{\lambda_1 + 1}{\lambda_1}\right) = \frac{\xi_m}{\frac{\partial \xi_m}{\partial \mathbf{L}_G}.\mathbf{L}_G + \frac{\partial \xi_m}{\partial \mathbf{D}_{p(t)}}.\mathbf{D}_{p(t)}}, \tag{3.14}$$

Hence, from Eqs. (3.12), (3.13), the final constitutive equations reduce to

$$\mathbf{T} = p\mathbf{I} + 2\rho \left[\left(\frac{\partial \hat{\psi}}{\partial I_{\mathbf{B}_{p(t)}}} + I_{\mathbf{B}_{p(t)}} \frac{\partial \hat{\psi}}{\partial II_{\mathbf{B}_{p(t)}}}\right) \mathbf{B}_{p(t)} - \frac{\partial \hat{\psi}}{\partial II_{\mathbf{B}_{p(t)}}} \mathbf{B}^2_{p(t)}\right] \\ + \left(\frac{\xi_m}{\frac{\partial \xi_m}{\partial \mathbf{L}_G}.\mathbf{L}_G + \frac{\partial \xi_m}{\partial \mathbf{D}_{p(t)}}.\mathbf{D}_{p(t)}}\right) \frac{\partial \xi_m}{\partial \mathbf{D}_{p(t)}}, \\ \mathbf{T} = \hat{\lambda}\mathbf{I} + \left(\frac{\xi_m}{\frac{\partial \xi_m}{\partial \mathbf{L}_G}.\mathbf{L}_G + \frac{\partial \xi_m}{\partial \mathbf{D}_{p(t)}}.\mathbf{D}_{p(t)}}\right) \mathbf{F}^{-T}_{\kappa_{p(t)}} \frac{\partial \xi_m}{\partial \mathbf{L}_G} \mathbf{F}^T_{\kappa_{p(t)}}, \tag{3.15}$$

where $\hat{\lambda} := \frac{\lambda_2}{\lambda_1}, p := \frac{\lambda_3}{\lambda_1}$ are Lagrange multipliers due to the constraint of incompressibility. At this point it appears that there are two constitutive relations – $(3.15)_a$, $(3.15)_b$ – for stress instead of just one! We would like to note that the two expressions – $(3.15)_a$, $(3.15)_b$ – have to be equated and simplified to obtain one expression for stress and an evolution equation for $\mathbf{B}_{p(t)}$. This will become clear within the context of the specific choices of $\psi$ and $\xi$ that are made in the next subsection. Furthermore, expression $(3.15)_b$ appears at first glance non-symmetric but is in fact symmetric for the specific form chosen for $\xi$ as we shall show in the next subsection.

3.2. **Specific model.** We now derive constitutive expression for the stress by choosing the stored energy to be

$$\hat{\psi}(I_{\mathbf{B}_{p(t)}}, II_{\mathbf{B}_{p(t)}}) = \frac{\mu}{2\rho}\left(I_{\mathbf{B}_{p(t)}} - 3\right), \tag{3.16}$$

and rate of dissipation to be of the form

$$\xi_m(\mathbf{B}_{p(t)}, \mathbf{D}_{p(t)}, \mathbf{D}_G) = \eta_p \mathbf{D}_{p(t)} . \mathbf{B}_{p(t)} \mathbf{D}_{p(t)} + \eta_G \mathbf{D}_G . \mathbf{B}_{p(t)} \mathbf{D}_G. \tag{3.17}$$

The stored energy chosen here is that for a neo-Hookean material with $\mu$ being its elastic modulus, whereas the rate of dissipation is similar to that of a "mixture" of two Newtonian-like fluids (in the sense that the dissipation is quadratic in the symmetric part of velocity gradient), whose dissipation also depends on the stretch (specifically the stretch from the natural configuration to the current configuration), with viscosities $\eta_G$ and $\eta_p$[4]. The former term on the right hand side of expression (3.17) is due to the dissipation during the motion from $\kappa_R$ to $\kappa_{p(t)}$ and the latter term is due to dissipation during the motion from $\kappa_{p(t)}$ to $\kappa_t$. Note that with the above choices for $\psi$ and $\xi$, as the body moves from the $\kappa_{p(t)}$ to $\kappa_t$, there is both storage (like a neo-Hookean solid) and dissipation (like a Newtonian-like fluid) of energy simultaneously and hence, $\kappa_{p(t)}$ evolves like the natural configuration of a Kelvin-Voigt-like solid (also see [20]) with respect to $\kappa_t$. Now, with this choice for $\psi$ and $\xi$, the constitutive relations given by Eq. (3.15) reduce to

$$\mathbf{T} = p\mathbf{I} + \mu \mathbf{B}_{p(t)} + \frac{\eta_p}{2} \left( \mathbf{B}_{p(t)} \mathbf{D}_{p(t)} + \mathbf{D}_{p(t)} \mathbf{B}_{p(t)} \right), \tag{3.18}$$

and

$$\mathbf{T} = \hat{\lambda} \mathbf{I} + \frac{\eta_G}{2} \mathbf{F}_{\kappa_{p(t)}}^{-T} \left( \mathbf{B}_{p(t)} \mathbf{D}_G + \mathbf{D}_G \mathbf{B}_{p(t)} \right) \mathbf{F}_{\kappa_{p(t)}}^{T}. \tag{3.19}$$

Now, from Eq. (3.18) and Eq. (3.19)

$$\left(p - \hat{\lambda}\right) \mathbf{I} + \mu \mathbf{B}_{p(t)} + \frac{\eta_p}{2} \left( \mathbf{B}_{p(t)} \mathbf{D}_{p(t)} + \mathbf{D}_{p(t)} \mathbf{B}_{p(t)} \right) =$$
$$\frac{\eta_G}{2} \mathbf{V}_{\kappa_{p(t)}}^{-1} \left( \mathbf{B}_{p(t)} \mathbf{D}_G + \mathbf{D}_G \mathbf{B}_{p(t)} \right) \mathbf{V}_{\kappa_{p(t)}}, \tag{3.20}$$

and using Eq. (2.15) and Eq. (3.5) in Eq. (3.20), we get

$$\left(p - \hat{\lambda}\right) \mathbf{I} + \mu \mathbf{B}_{p(t)} + \frac{\eta_p}{4} \left( \mathbf{V}_{\kappa_{p(t)}} \dot{\mathbf{B}}_{p(t)} \mathbf{V}_{\kappa_{p(t)}}^{-1} + \mathbf{V}_{\kappa_{p(t)}}^{-1} \dot{\mathbf{B}}_{p(t)} \mathbf{V}_{\kappa_{p(t)}} \right) =$$
$$\frac{\eta_G}{4} \mathbf{V}_{\kappa_{p(t)}}^{-1} \left[ -\mathbf{V}_{\kappa_{p(t)}} \overset{\nabla}{\mathbf{B}}_{p(t)} \mathbf{V}_{\kappa_{p(t)}}^{-1} - \mathbf{V}_{\kappa_{p(t)}}^{-1} \overset{\nabla}{\mathbf{B}}_{p(t)} \mathbf{V}_{\kappa_{p(t)}} \right] \mathbf{V}_{\kappa_{p(t)}}. \tag{3.21}$$

Pre-multiplying and post-multiplying Eq. (3.21) by $\mathbf{V}_{\kappa_{p(t)}}$, we have

$$\left(p - \hat{\lambda}\right) \mathbf{B}_{p(t)} + \mu \mathbf{B}_{p(t)}^2 + \frac{\eta_p}{4} \left( \mathbf{B}_{p(t)} \dot{\mathbf{B}}_{p(t)} + \dot{\mathbf{B}}_{p(t)} \mathbf{B}_{p(t)} \right) =$$
$$\frac{\eta_G}{4} \mathbf{V}_{\kappa_{p(t)}}^{-1} \left( -\mathbf{B}_{p(t)} \overset{\nabla}{\mathbf{B}}_{p(t)} - \overset{\nabla}{\mathbf{B}}_{p(t)} \mathbf{B}_{p(t)} \right) \mathbf{V}_{\kappa_{p(t)}}. \tag{3.22}$$

Also, from Eq. (3.21)

$$\left(p - \hat{\lambda}\right) = -\frac{1}{3} \left[ \frac{\eta_G}{2} tr\left(\overset{\nabla}{\mathbf{B}}_{p(t)}\right) + \frac{\eta_p}{2} tr\left(\dot{\mathbf{B}}_{p(t)}\right) + \mu tr\left(\mathbf{B}_{p(t)}\right) \right]. \tag{3.23}$$

Eq. (3.19) can be re-written as

$$\mathbf{T} = \hat{\lambda} \mathbf{I} + \frac{\eta_G}{4} \mathbf{V}_{\kappa_{p(t)}}^{-1} \left[ -\mathbf{V}_{\kappa_{p(t)}} \overset{\nabla}{\mathbf{B}}_{p(t)} \mathbf{V}_{\kappa_{p(t)}}^{-1} - \mathbf{V}_{\kappa_{p(t)}}^{-1} \overset{\nabla}{\mathbf{B}}_{p(t)} \mathbf{V}_{\kappa_{p(t)}} \right] \mathbf{V}_{\kappa_{p(t)}}. \tag{3.24}$$

---

[4]Of course, one can choose the rate of dissipation to be quadratic in the symmetric part of velocity gradient without any stretch dependence, for example, $\xi = \eta_p \mathbf{D}_{p(t)} . \mathbf{D}_{p(t)} + \eta_G \mathbf{D}_G . \mathbf{D}_G$. The resulting model would be a variation of the current model. See Section 7 for the calculations using such a rate of dissipation.

Eq. (3.22) and Eq. (3.24) are the final constitutive relations with Eq. (3.22) togther with Eq. (3.23) being the evolution equation for the natural configuration $\kappa_{p(t)}$.

In what follows, we shall also show that our model can be reduced to both the Maxwell-like fluid model and Kelvin-Voigt-like solid model under certain assumptions for the material parameters that are involved. We shall solve the problems of creep and stress relaxation by considering homogeneous uniaxial extension. Based on the results for creep and stress relaxation, and following the definitions given in [31, 32] for a fluid-like body and a solid-like body, we shall show that our model is a fluid-like model when none of the material parameters are ignored.

3.3. **Limiting cases.** By setting $\eta_p$ to zero, one can see from Eq. (3.20) that $\mathbf{B}_{p(t)}$ and $\mathbf{D}_G$ have the same eigen-vectors and hence commute. Using this fact, Eq. (3.18) and Eq. (3.22) reduce to the Maxwell-like fluid model developed by Rajagopal and Srinivasa (see Eqs. 40–42 in [5]). On the other hand, if we assume $\eta_G$ goes to infinity, with $\eta_p, \mu$ and all other kinematical quantities remaining finite, then as the deviatoric part of Eq. (3.19) is finite, from Eqs. (3.19), (3.10), we must have $\mathbf{D}_G \to 0$. This also implies that $\mathbf{G}$ is pure rotation and hence $\mathbf{B}_{p(t)} \to \mathbf{B}_{\kappa_R}$, $\mathbf{D}_{p(t)} \to \mathbf{D}$. Thus, Eq. (3.18) reduces to

$$\mathbf{T} = p\mathbf{I} + \mu \mathbf{B}_{\kappa_R} + \frac{\eta_p}{2} \left( \mathbf{B}_{\kappa_R} \mathbf{D} + \mathbf{D} \mathbf{B}_{\kappa_R} \right). \tag{3.25}$$

This is a generalized Kelvin-Voigt solid model.

## 4. Application of the model

4.1. **Creep.** In this subsection we shall solve the problem of homogeneous extension under constant stress and show that our model is a fluid-like model when none of the material moduli are ignored. Before we go into the details of the problem, we would like to note that for steady problems wherein the deformation is homogeneous, the results for our model would be same as that for a Maxwell-like fluid with stretch dependent relaxation as formulated in [5]. This is because $\dot{\mathbf{B}}_{p(t)} = 0$ (since $\frac{\partial \mathbf{B}_{p(t)}}{\partial t} = 0$, due to the assumption that the deformation is steady and $grad(\mathbf{B}_{p(t)}) = 0$ as the deformation is homogeneous) and the constitutive relations in Eqs. (3.22), (3.24) reduce to Eqs. (36), (39) in [5].

Now, in the case of time dependent homogeneous extension:

$$x = \lambda(t)X, \quad y = \frac{1}{\sqrt{\lambda(t)}}Y, \quad z = \frac{1}{\sqrt{\lambda(t)}}Z, \tag{4.1}$$

the deformation gradient with respect to the reference configuration is given by

$$\mathbf{F}_{\kappa_R} = \text{diag}\left\{ \lambda(t), \frac{1}{\sqrt{\lambda(t)}}, \frac{1}{\sqrt{\lambda(t)}} \right\}. \tag{4.2}$$

Hence, the velocity gradient is given by

$$\mathbf{L} = \text{diag}\left\{ \frac{\dot{\lambda}}{\lambda}, -\frac{\dot{\lambda}}{2\lambda}, -\frac{\dot{\lambda}}{2\lambda} \right\}. \tag{4.3}$$

Now, we shall assume that

$$\mathbf{B}_{p(t)} = \text{diag}\left\{ B(t), \frac{1}{\sqrt{B(t)}}, \frac{1}{\sqrt{B(t)}} \right\}, \tag{4.4}$$

and hence
$$\dot{\mathbf{B}}_{p(t)} = \text{diag}\left\{\dot{B}(t), -\frac{\dot{B}(t)}{2B^{3/2}(t)}, -\frac{\dot{B}(t)}{2B^{3/2}(t)}\right\}, \quad (4.5)$$

$$\overset{\nabla}{\mathbf{B}}_{p(t)} = \text{diag}\left\{\dot{B}(t) - \frac{2B(t)\dot{\lambda}(t)}{\lambda(t)}, -\frac{\dot{B}(t)}{2B^{3/2}(t)} + \frac{\dot{\lambda}(t)}{\sqrt{B(t)}\lambda(t)}, -\frac{\dot{B}(t)}{2B^{3/2}(t)} + \frac{\dot{\lambda}(t)}{\sqrt{B(t)}\lambda(t)}\right\}, \quad (4.6)$$

and
$$\mathbf{V}_{\kappa_{p(t)}} = \text{diag}\left\{\sqrt{B}(t), \frac{1}{B^{1/4}(t)}, \frac{1}{B^{1/4}(t)}\right\}. \quad (4.7)$$

For the case of homogeneous extension, Eq. (3.18) and Eq. (3.19) reduce to
$$\mathbf{T} = p\mathbf{I} + \mu\mathbf{B}_{p(t)} + \eta_p \mathbf{B}_{p(t)}\mathbf{D}_{p(t)},$$
$$\mathbf{T} = \hat{\lambda}\mathbf{I} + \eta_G \mathbf{B}_{p(t)}\mathbf{D}_G, \quad (4.8)$$

and the final constitutive relations Eq. (3.22) and Eq. (3.24) become
$$\mathbf{T} = \hat{\lambda}\mathbf{I} - \frac{\eta_G}{2}\mathbf{V}_{\kappa_{p(t)}} \overset{\nabla}{\mathbf{B}}_{p(t)} \mathbf{V}_{\kappa_{p(t)}}^{-1}, \quad (4.9)$$

$$(p - \hat{\lambda})\mathbf{I} + \mu\mathbf{B}_{p(t)} = -\frac{\eta_G}{2}\overset{\nabla}{\mathbf{B}}_{p(t)} - \frac{\eta_p}{2}\dot{\mathbf{B}}_{p(t)}, \quad (4.10)$$

where
$$(p - \hat{\lambda}) = -\frac{3\mu}{tr(\mathbf{B}_{p(t)}^{-1})}. \quad (4.11)$$

On substituting Eqs. (4.5), (4.6) into Eq. (4.10), we arrive at
$$\dot{B}(t) = \frac{2\eta_G}{\eta_G + \eta_p}\frac{B(t)\dot{\lambda}(t)}{\lambda(t)} - \frac{4\mu}{\eta_G + \eta_p}\left[\frac{B^{5/2}(t) - B(t)}{1 + 2B^{3/2}(t)}\right]. \quad (4.12)$$

Using Eqs. (4.6), (4.7) in Eq. (4.9) and using the fact that lateral surfaces are traction free, we find that
$$T_{11} = \frac{\eta_G}{2}\left(1 + \frac{1}{2B^{3/2}(t)}\right)\left(\frac{2B(t)\dot{\lambda}(t)}{\lambda(t)} - \dot{B}(t)\right). \quad (4.13)$$

Solving Eqs. (4.12), (4.13) simultaneously, we obtain
$$\dot{B}(t) = \frac{2T_{11}}{\eta_p\left(1 + \frac{1}{2B^{3/2}(t)}\right)} - \frac{4\mu}{\eta_p}\left(\frac{B^{5/2}(t) - B(t)}{1 + 2B^{3/2}(t)}\right), \quad (4.14)$$

$$\frac{\dot{\lambda}(t)}{\lambda(t)} = \frac{2(\eta_p + \eta_G)}{\eta_G \eta_p}\left(\frac{T_{11}B^{1/2}(t)}{1 + 2B^{3/2}(t)}\right) - \frac{2\mu}{\eta_p}\left(\frac{B^{3/2}(t) - 1}{1 + 2B^{3/2}(t)}\right), \quad (4.15)$$

which can be re-written as
$$\frac{dB}{d\bar{t}} = \frac{2\bar{T}_{11}}{\left(1 + \frac{1}{2B^{3/2}(t)}\right)} - 4\left(\frac{B^{5/2}(t) - B(t)}{1 + 2B^{3/2}(t)}\right), \quad (4.16)$$

$$\frac{1}{\lambda}\frac{d\lambda}{d\bar{t}} = 2\left(\frac{1}{\bar{\eta}} + 1\right)\left(\frac{\bar{T}_{11}B^{1/2}(t)}{1 + 2B^{3/2}(t)}\right) - 2\left(\frac{B^{3/2}(t) - 1}{1 + 2B^{3/2}(t)}\right), \quad (4.17)$$

where $\bar{t} = \frac{t\mu}{\eta_p}$ is a non-dimensional time, $\bar{T}_{11} = \frac{T_{11}}{\mu}$ is a non-dimensional stress, and $\bar{\eta} = \frac{\eta_p}{\eta_G}$ is the ratio of viscosities $(\eta_p, \eta_G)$. With known values for $\bar{T}_{11}$ and $\bar{\eta}$, the ODEs – Eqs. (4.16), (4.17) – were solved simultaneously using "ode45" solver in Matlab. The initial conditions $B(0) = \lambda(0) = 1$ were used for the loading process. In the case of unloading, the initial values were set to the values of $B$, $\lambda$ evaluated at the end of loading process. Before we discuss our results, the reader should note that $B$ is the square of the stretch from the natural configuration to the current configuration and $\lambda$ is the stretch from the reference configuration to the current configuration or the total stretch.

We shall now discuss the numerical results obtained for the problem of creep. Fig. (2) portrays the results (specifically, $B$ as a function of $\bar{t}$ and $\lambda$ as function of $\bar{t}$) for $\bar{T}_{11} = 1$ for various values of $\bar{\eta}$. For all the cases of $\bar{\eta}$, the evolution of the natural configuration with respect to the current configuration (or $B$ as a function of $\bar{t}$) is the same. This can also be seen from Eqs. (4.16),(4.17) as the evolution of $B$ in time depends only on the value of $\bar{T}_{11}$, and does not depend on $\bar{\eta}$. Also, as seen in Fig. (2)(a), on removal of the load, the stretch in the body does not return to unity but exhibits a permanent residual stretch, which is the typical behavior of a fluid-like body (see [31, 32]). In addition, from Fig. (2) and Fig. (3), by increasing $\bar{T}_{11}$ (from 1 to 5) for fixed $\bar{\eta}$ (here the value is 10) – the maximum value of $B$ increases, the maximum value for $\lambda$ increases and the permanent residual stretch also increases.

Now, as discussed earlier in section 3.2, the natural configuration evolves like a Kelvin-Voigt-like solid with respect to the current configuration. Fig. (2)(b) reiterates this fact since $B$ varies with $t$ in a fashion similar to the stretch as a function of time for a Kelvin-Voigt solid in a creep experiment. Fig. (2)(a) also shows that, as $\bar{\eta}$ increases, the maximum value for total stretch ($\lambda$) decreases. This shows that for a fixed amount of loading time (for instance for $\bar{t} = 10$) as the value for $\bar{\eta}$ increases, the value for maximum total stretch decreases and so the rate of relaxation decreases. As $\bar{\eta}$ increases one would expect that our model would tend to a solid model in the limit of $\bar{\eta} \to \infty$. To see this, in our numerical simulations, we set $\bar{\eta}$ to a very large number (specifically $\bar{\eta} = 100000$), with $\bar{T}_{11} = 1$, and we found $\sqrt{B}$ and $\lambda$ to follow the same trend in time (see Fig. 4). This is an extreme case when the natural configuration and the reference configuration tend to being the same, reducing our model to Kelvin-Voigt-like solid model.

4.2. **Stress under constant strain rate.** If we define $\epsilon := ln\lambda$ as our strain measure, then Eq. (4.12) reduces to

$$\dot{B}(t) = \frac{2\eta_G}{\eta_G + \eta_p} B(t)\dot{\epsilon} - \frac{4\mu}{\eta_G + \eta_p}\left[\frac{B^{5/2}(t) - B(t)}{1 + 2B^{3/2}(t)}\right]. \tag{4.18}$$

Using Eq. (4.18) in Eq. (4.13), we get

$$T_{11} = \eta_G\left(1 + \frac{1}{2B^{3/2}(t)}\right)\left[\frac{\eta_p}{\eta_p + \eta_G}B(t)\dot{\epsilon} + \frac{2\mu}{\eta_G + \eta_p}\frac{B^{5/2}(t) - B(t)}{1 + 2B^{3/2}(t)}\right]. \tag{4.19}$$

Upon non-dimensionalizing Eqs. (4.18), (4.19), we arrive at

$$\frac{dB}{d\bar{t}} = \frac{2\bar{\eta}B}{\bar{\eta} + 1}\frac{d\epsilon}{d\bar{t}} - \frac{4}{\bar{\eta} + 1}\left(\frac{B^{5/2} - B}{1 + 2B^{3/2}}\right), \tag{4.20}$$

$$\bar{T}_{11} = \left(1 + \frac{1}{2B^{3/2}}\right)\left[\frac{B\bar{\eta}}{1 + \bar{\eta}}\frac{d\epsilon}{d\bar{t}} + \frac{2}{1 + \bar{\eta}}\left(\frac{B^{5/2} - B}{1 + 2B^{3/2}}\right)\right]. \tag{4.21}$$

With known values for $\frac{d\epsilon}{d\bar{t}}$, $B$ can be evaluated using Eq. (4.20) and then from Eq. (4.21), $\bar{T}_{11}$ can be calculated. For stress relaxation, we set $\frac{d\epsilon}{d\bar{t}}$ to zero and solved the ODEs – Eqs. (4.20), (4.21) – numerically using the solver "ode45" in Matlab. Fig. (5) displays the plots for $B$, and $\bar{T}_{11}$ as functions of $\bar{t}$ for various values of $\bar{\eta}$. The initial condition for B was chosen to be 2. Fig. (5)(b) also shows that the stress finally relaxes to zero. This is a characteristic of a fluid-like material. In addition, one can also see from Fig. (5) that as the ratio of viscosities $\bar{\eta}$ increases, the time of relaxation increases. Also from Eq. (4.21), as $\bar{\eta} \to \infty$, $\bar{T}_{11} \to 0$ and hence cannot stress relax. This is because as $\bar{\eta} \to \infty$, the model behaves like a Kelvin-Voigt-like solid model as can be gleaned from the creep response discussed in Section (4.1). In addition, one can also notice from Fig. (5)(b) that the initial value for $\bar{T}_{11}$ decreases as $\bar{\eta}$ increases.

## 5. Model reduction in one dimension

If the elastic strain is small in the sense that

$$\|\mathbf{B}_{p(t)} - \mathbf{I}\| = O(\epsilon), \quad \epsilon \ll 1, \tag{5.1}$$

then Eq. (4.11) reduces to (see [5])

$$\left(p - \hat{\lambda}\right) = -\mu, \tag{5.2}$$

and hence Eq. (4.8) reduce to

$$\begin{aligned} \mathbf{T} &= \hat{\lambda}\mathbf{I} + \eta_G \mathbf{B}_{p(t)} \mathbf{D}_G, \\ \mu\left(\mathbf{B}_{p(t)} - \mathbf{I}\right) &= \mathbf{B}_{p(t)} \left(\eta_G \mathbf{D}_G - \eta_p \mathbf{D}_{p(t)}\right). \end{aligned} \tag{5.3}$$

If $\lambda_i$ $(i = G, p)$ is the one-dimensional stretch and $\epsilon = ln\lambda_i$ is the true strain, then in one dimension, Eq. (5.3) reduces to

$$\begin{aligned} \sigma &= \hat{\eta}_G \frac{\dot{\lambda}_G}{\lambda_G}, \\ \mu\left(\lambda_p^2 - 1\right) &= \hat{\eta}_G \frac{\dot{\lambda}_G}{\lambda_G} - \hat{\eta}_p \frac{\dot{\lambda}_p}{\lambda_p}, \end{aligned} \tag{5.4}$$

where $\hat{\eta}_i = \eta_i \lambda_p^2$ $(i = G, p)$ are the stretch dependent viscosities and $\sigma$ is the one dimensional stress. Eq. (5.4) under the assumption that $\epsilon_i \ll 1$ $(i = G, p)$ reduces to

$$\sigma = \hat{\eta}_G \dot{\epsilon}_G, \quad 2\mu\epsilon_p = \hat{\eta}_G \dot{\epsilon}_G - \hat{\eta}_p \dot{\epsilon}_p. \tag{5.5}$$

Eqn. (5.5) can also be obtained if we have a dashpot with $\hat{\eta}_G$ as the viscosity and a Kelvin-Voigt element, with a spring constant $2\mu$ and viscosity $\hat{\eta}_p$, in series (see Fig. 1). Thus, our model in Section 3.2 reduces to a dashpot and a Kelvin-Voigt element in series whose viscosities are stretch dependent.

## 6. Concluding Remarks

For the model developed in this paper, we solved the problems of creep and stress relaxation and based on the results showed that our model is a fluid-like model. We have also shown that under certain idealizations, our model reduces to both a Maxwell-like fluid and a Kelvin-Voigt-like solid. We have also shown that our model is a three dimesional generalization of the one dimensional model with a dashpot and Kelvin-Voigt element in series. We would also like to note that the general framework developed in Section 3.1 can be extended to model visco-elasto-plastic response by choosing the rate of dissipation for the motion from $\kappa_R$ to $\kappa_{p(t)}$ to be of the form given in Eq. (39) or Eq. (40) in [9]. For example, one can choose the rate of dissipation to be of the form

$$\xi_m(\mathbf{B}_{p(t)}, \mathbf{D}_{p(t)}, \mathbf{D}_G) = Y\left(\mathbf{D}_{p(t)}.\mathbf{B}_{p(t)}\mathbf{D}_{p(t)}\right)^{1/2} + \eta_G \mathbf{D}_G.\mathbf{B}_{p(t)}\mathbf{D}_G, \tag{6.1}$$

where $Y$ is a material constant, along with the stored energy in Eq. (3.16), and then maximize the rate of dissipation under appropriate constraints.

## 7. Appendix

In this section, we shall derive the constitutive relations when the rate of dissipation is of the form

$$\xi_m(\mathbf{D}_{p(t)}, \mathbf{D}_G) = \eta_p \mathbf{D}_{p(t)}.\mathbf{D}_{p(t)} + \eta_G \mathbf{D}_G.\mathbf{D}_G. \tag{A.1}$$

with the stored energy given by Eq. (3.16). Now, Eq. (3.15) reduces to

$$\mathbf{T} = p\mathbf{I} + \mu \mathbf{B}_{p(t)} + \eta_p \mathbf{D}_{p(t)}, \tag{A.2}$$

and

$$\mathbf{T} = \hat{\lambda}\mathbf{I} + \eta_G \mathbf{V}_{\kappa_{p(t)}}^{-1} \mathbf{D}_G \mathbf{V}_{\kappa_{p(t)}}. \tag{A.3}$$

Also, from Eqs. (A.2) and (A.3),

$$(p - \hat{\lambda}) = -\frac{\mu}{3} tr(\mathbf{B}_{p(t)}). \tag{A.5}$$

Equating Eqs. (A.2) and (A.3), pre-multiplying and post-multiplying by $\mathbf{V}_{\kappa_{p(t)}}$, and using Eqs. (2.15), (3.5), we get

$$(p - \hat{\lambda})\mathbf{B}_{p(t)} + \mu \mathbf{B}_{p(t)}^2 + \frac{\eta_p}{2}\dot{\mathbf{B}}_{p(t)} = -\frac{\eta_G}{2}\mathbf{V}_{\kappa_{p(t)}}^{-1} \overset{\nabla}{\mathbf{B}}_{p(t)} \mathbf{V}_{\kappa_{p(t)}}. \tag{A.6}$$

It follows that, the constitutive relations for the choices of Eq. (3.16) and Eq. (A.1) for the specific Helmholtz potential and the rate of dissipation, are

$$\mathbf{T} = p\mathbf{I} + \mu \mathbf{B}_{p(t)} + \frac{\eta_p}{2}\mathbf{V}_{\kappa_{p(t)}}^{-1} \dot{\mathbf{B}}_{p(t)} \mathbf{V}_{\kappa_{p(t)}}^{-1}, \tag{A.7}$$

and

$$-\frac{\eta_G}{2}\mathbf{V}_{\kappa_{p(t)}}^{-1} \overset{\nabla}{\mathbf{B}}_{p(t)} \mathbf{V}_{\kappa_{p(t)}} = -\frac{\mu}{3} tr(\mathbf{B}_{p(t)})\mathbf{B}_{p(t)} + \mu \mathbf{B}_{p(t)}^2 + \frac{\eta_p}{2}\dot{\mathbf{B}}_{p(t)}. \tag{A.8}$$

Eq. (A.8) is the evolution equation of the natural configuration.

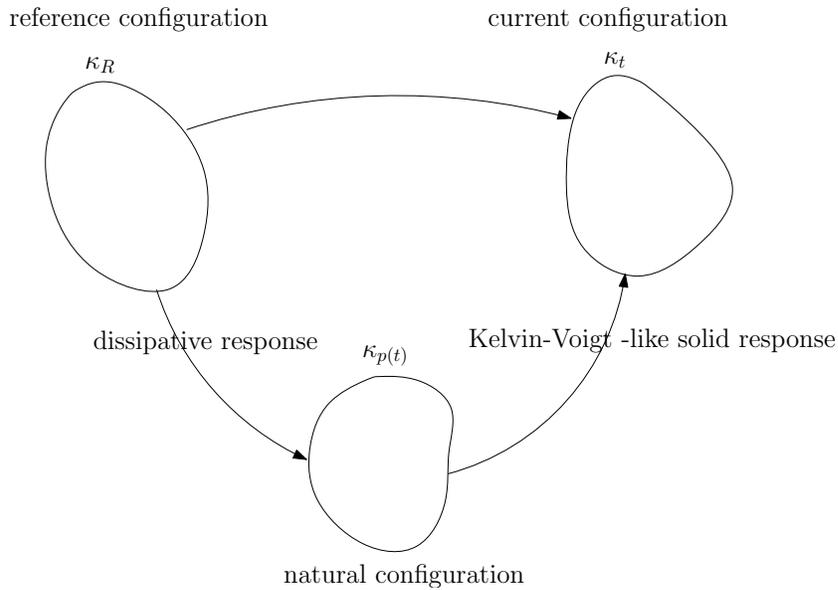

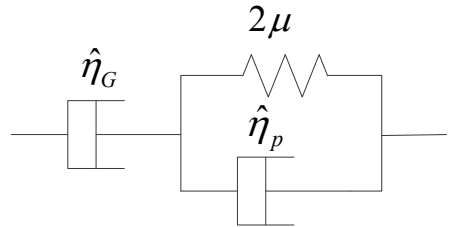

FIGURE 1. Schematic of the natural configuration $\kappa_{p(t)}$ corresponding to the current configuration $\kappa_t$ and the relevant mappings from the tangent spaces of the same material point in $\kappa_R$, $\kappa_t$ and $\kappa_{p(t)}$ (above). The response from the natural configuration $\kappa_{p(t)}$ is like a Kelvin-Voigt solid and the response of $\kappa_{p(t)}$ from the reference configuration $\kappa_R$ is purely dissipative. The corresponding one-dimesional spring dashpot analogy where a dashpot is in series with a Kelvin-Voigt element (below).

Satish Karra, Texas A&M University, Department of Mechanical Engineering, 3123 TAMU, College Station TX 77843-3123, United States of America
  *E-mail address*: `satkarra@tamu.edu`

K. R. Rajagopal (Corresponding author), Texas A&M University, Department of Mechanical Engineering, 3123 TAMU, College Station TX 77843-3123, United States of America, Phone: 1-979-862-4552, Fax: 1-979-845-3081
  *E-mail address*: `krajagopal@tamu.edu`

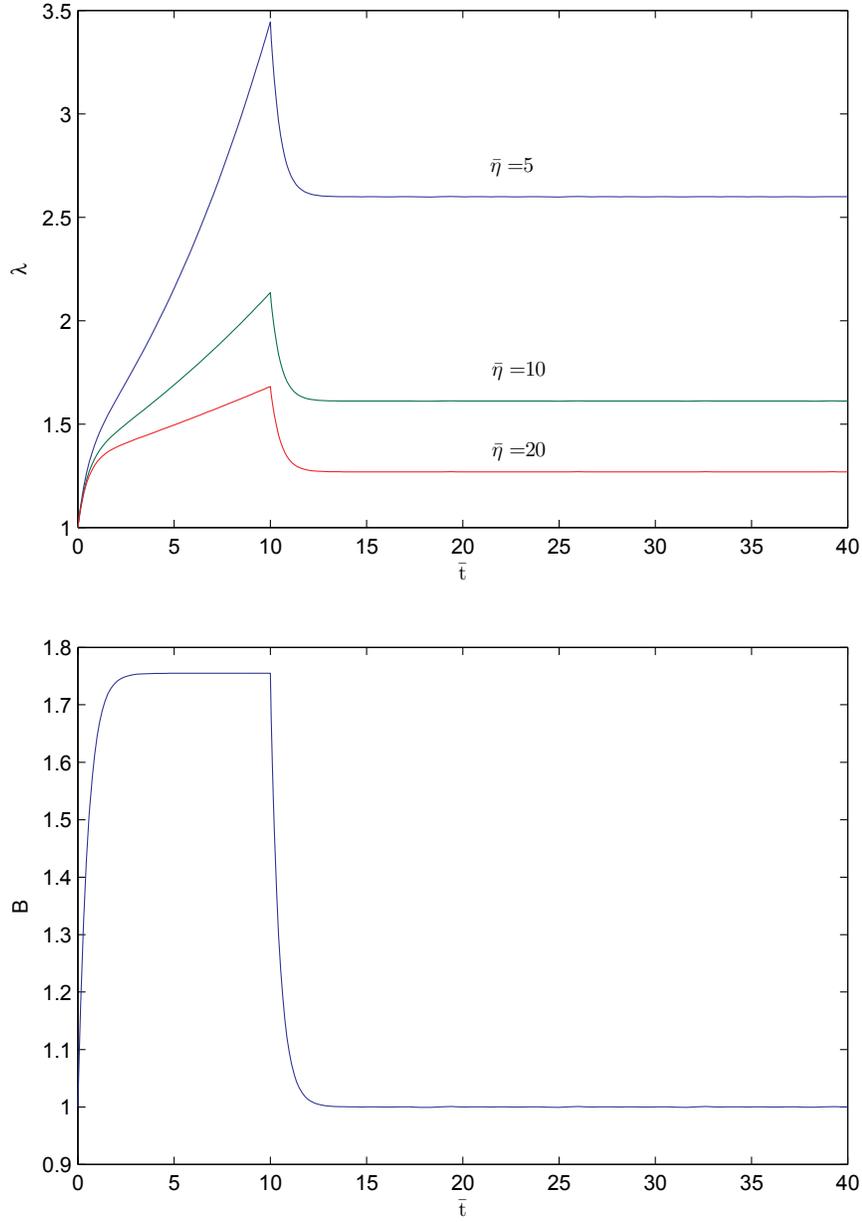

FIGURE 2. Overall stretch of the current configuration from the reference configuration ($\lambda$) and square of the stretch from the natural configuration to the current configuration ($B$) as a function of non-dimensional time $\bar{t}$ for the creep experiment. For the loading process, $\bar{T}_{11} = 1$ and the unloading starts at $t = 10$. Plots for $\bar{\eta} = 5,\ 10,\ 20$ are shown.

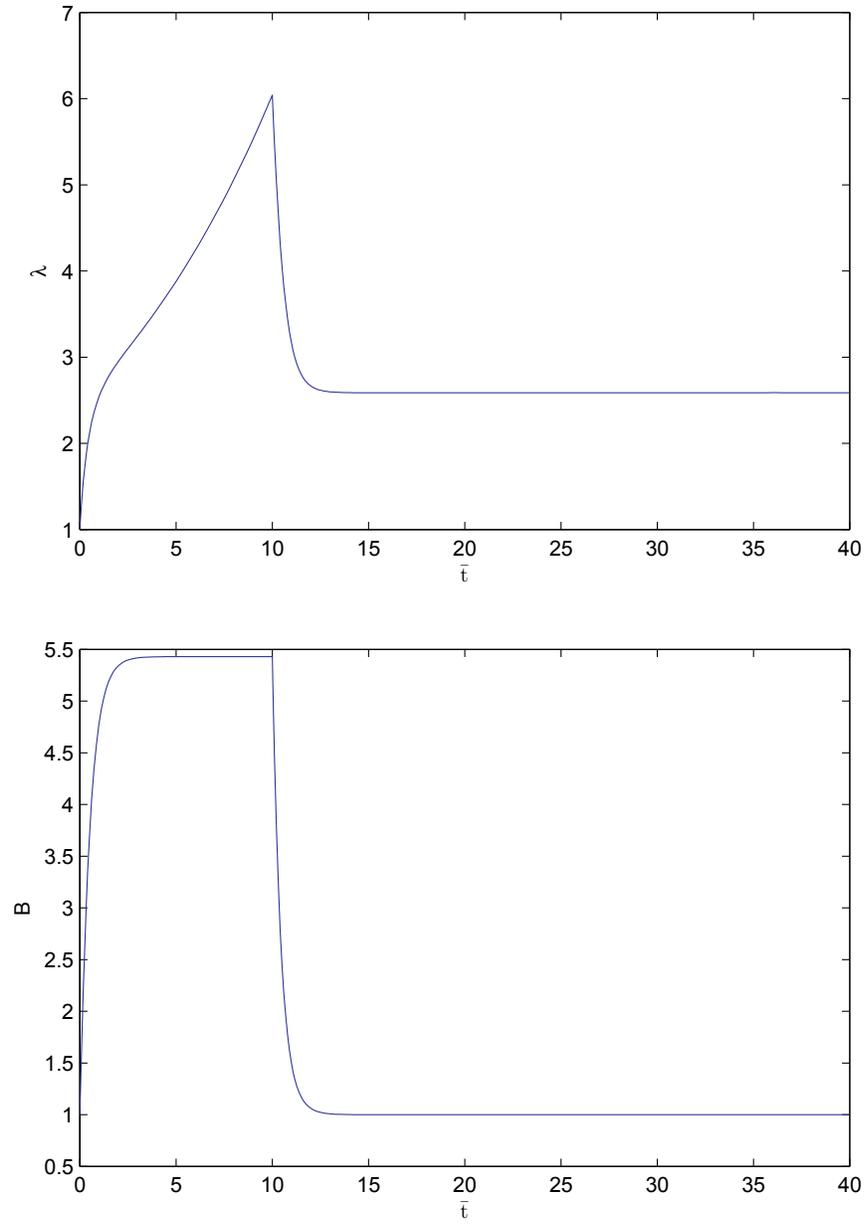

FIGURE 3. Overall stretch of the current configuration from the reference configuration ($\lambda$) and square of the stretch from the natural configuration to the current configuration ($B$) as a function of non-dimensional time $\bar{t}$ for the creep experiment. For this case, the non-dimensional stress $\bar{T}_{11} = 5$ and the ratio of the viscosities $\bar{\eta} = 10$. The unloading for the creep experiment starts at $t = 10$.

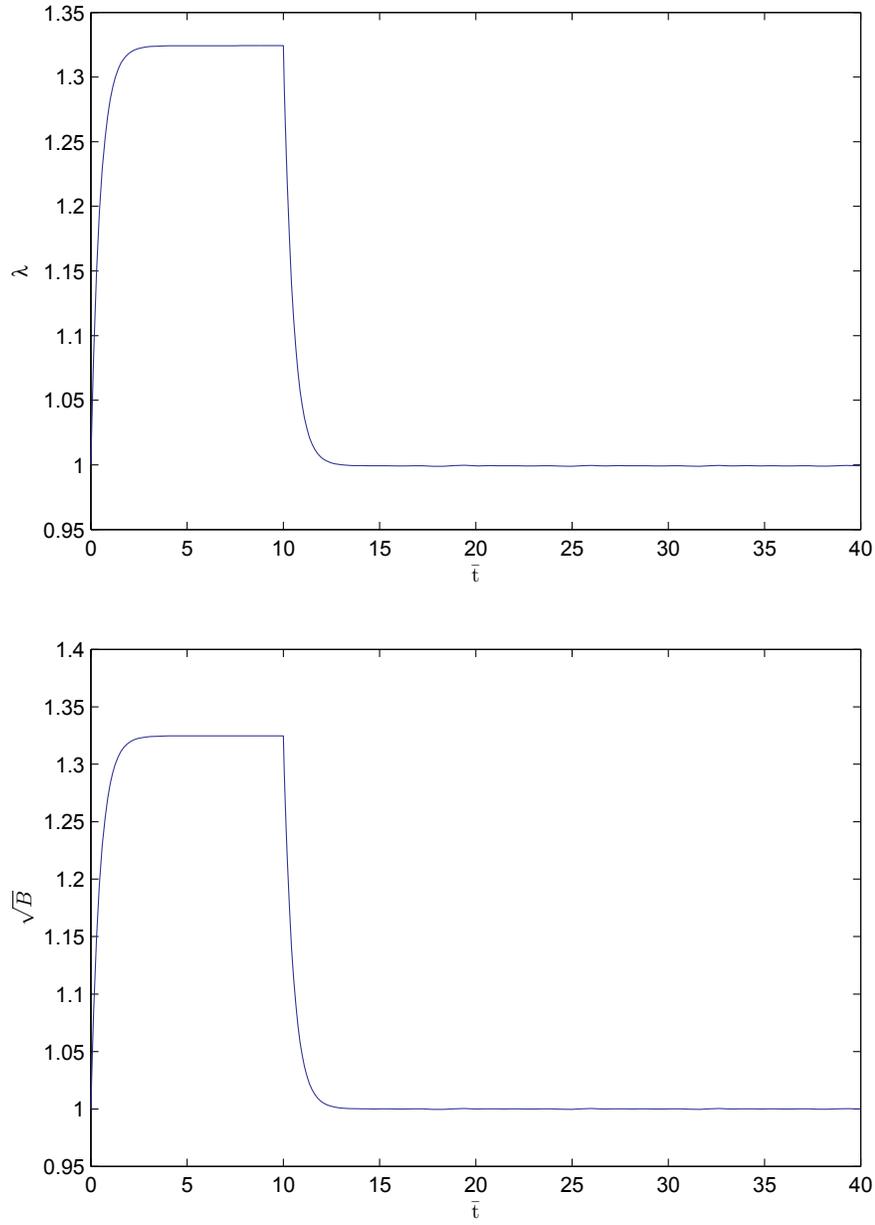

FIGURE 4. Overall stretch of the current configuration from the reference configuration ($\lambda$) and stretch from the natural configuration to the current configuration ($\sqrt{B}$) as a function of non-dimensional time $\bar{t}$ for the creep experiment. For this case, the non-dimensional stress $\bar{T}_{11} = 1$ and the ratio of the viscosities $\bar{\eta} = 100000$. Unloading for the creep experiment starts at $t = 10$.

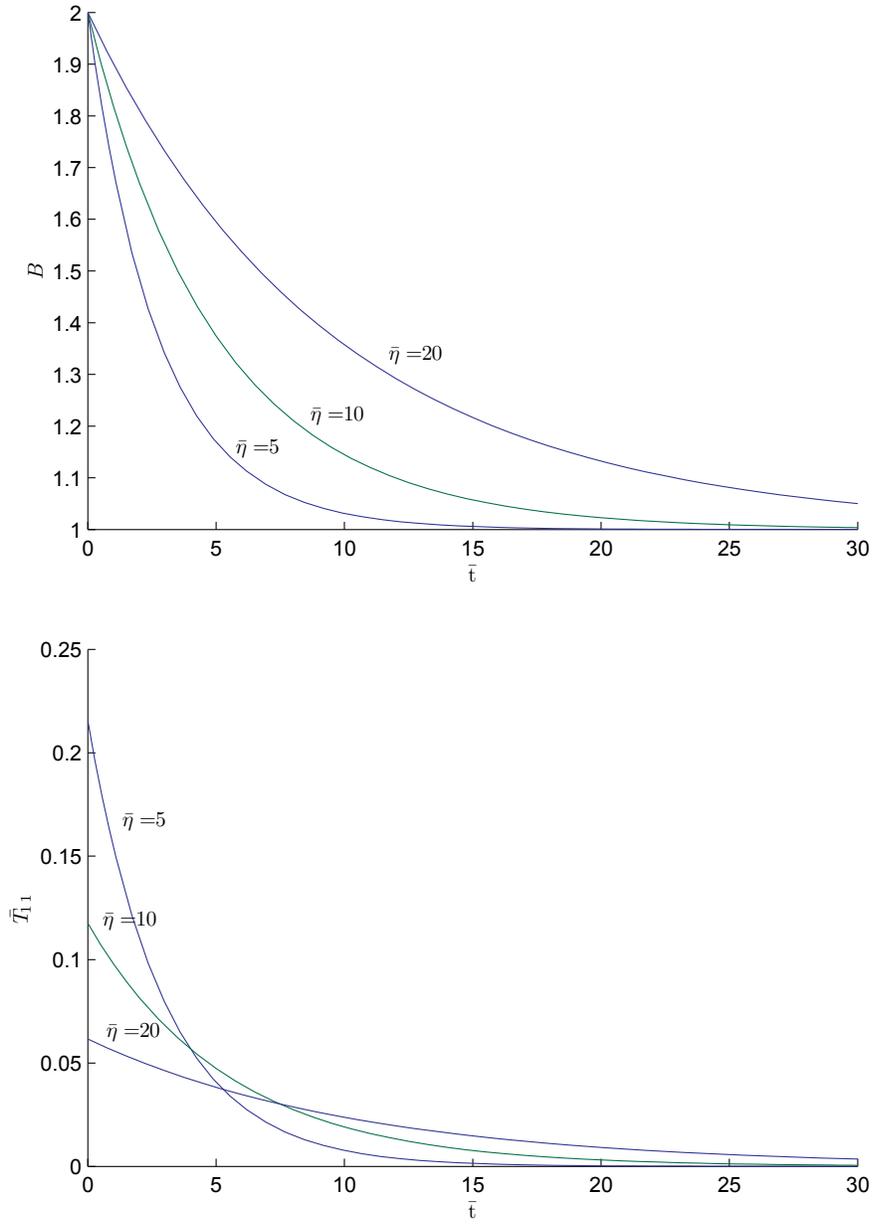

FIGURE 5. Square of the stretch from the natural configuration to the current configuration ($B$), non-dimensional stress ($\bar{T}_{11}$) plotted as functions of non-dimensional time $\bar{t}$ for various values of $\bar{\eta}$, in the stress relaxation experiment. The initial condition for B was chosen as 2.